\documentclass[11pt,draft,twoside]{article}
\usepackage[cp1250]{inputenc}
\usepackage{latexsym,amssymb,amsmath}
\title{About the embedding of Moufang loops in alternative algebras II}
\author{Sandu N. I.}
\date{}
\begin{document}
\maketitle

\begin{abstract}
It is known that with precision till isomorphism that only and
only loops $M(F) = M_0(F)/<-1>$, where $M_0(F)$ denotes the loop,
consisting from elements of all matrix Cayley-Dickson algebra
$C(F)$ with norm $1$, and $F$ be a subfield of arbitrary fixed
algebraically closed field, are simple non-associative Moufang
loops. In this paper it is proved that the simple loops $M(F)$
they and only they are not embedded into a loops of invertible
elements of any unitaly alternative algebras if $\text{char} F
\neq 2$ and $F$ is closed under square root operation. For the
remaining Moufang loops such an embedding is possible. Using this
embedding it is quite simple to prove the well-known finding: the
finite Moufang $p$-loop is centrally nilpotent.
\smallskip\\

\smallskip\\
\textit{Classification:} 17D05, 20N05.

\end{abstract}
\bigskip

\section{Introduction}

For an alternative algebra $A$ with the unit $1$ the set $U(A)$ of
all invertible elements of $A$ forms a Moufang loop with respect
to multiplication [1]. In [2] it is proved that any relative free
Moufang loop can be embedded into a loop of type $U(A)$. This
gives a positive answer to the question, raised by I. P. Shestakov
in [3]. But in [4] E. G. Goodaire raises a broader question: is it
true that any Moufang loop can be embedded into a loop of type
$U(A)$ for a suitable unital alternative algebra $A$? Analogical
question for commutative Moufang loops is raised by A. N. Grishkov
in the works of Loops'03 Conference (Prague, 2003). A positive
answer to Goodaire's question was announced in [5]. Here the
answer to this question is negative: in [3] there is constructed
an example of Moufang loops which are not embedded into loops of
type $U(A)$. We mentioned that by Theorem from [6] these Moufang
loops are simple.

It is known [7, 8] that with precision till isomorphism that only
and only loops $M(F) = M_0(F)/<-1>$, where $M_0(F)$ denotes the
loop, consisting from elements of all matrix Cayley-Dickson
algebra $C(F)$ with norm $1$, and $F$ be a subfield of arbitrary
fixed algebraically closed field, are simple non-associative
Moufang loops. In this paper, which is a continuation of [2],
Goodaire's question is completely settled, namely. It is proved
that in class of all  Moufang loops the simple loops $M(F)$ they
and only they are not embedded into a loops of type $U(A)$ of any
unitaly alternative algebras $A$ if $\text{char} F \neq 2$ and $F$
is closed under square root operation. The remaining loops can be
embedded into  loops of type $U(A)$. In particular, any non-simple
loops $Q$ can be embedded into a loop $U(FQ)$, where $F$ is an
arbitrary field and $FQ$ is the $''$loop algebra$''$ of $Q$.

In [9, 10], using the strong apparatus of finite group theory, and
quite unwieldy it is proved that the finite Moufang $p$-loop is
centrally nilpotent. In our work, a quite simple proof  of this
result is stated using the embedding of Moufang loops in
alternative algebras. Besides, in the proof only such a fact is
used:  the alternative $F$-algebra, generated  as $F$-module by
nilpotent generators  is nilpotent.

For basic definitions and properties of loops, see [11, 12], of
alternative algebras, see [13], and of fields, see [14].

\section{Centrally nilpotent loops}

Loop $(Q,\cdot) \equiv Q$ is called \textit{$IP$-loop} if the laws
$^{-1}x\cdot xy = yx\cdot x^{-1} = y$ are  true  in  it,  where
$^{-1}xx = xx^{-1} = 1$. In $IP$-loops $^{-1}x = x^{-1}$  and
$(xy)^{-1} = y^{-1}x^{-1}$. The \textit{Moufang loop} is defined
by identity
$$(xy\cdot x)z = x(y\cdot xz). \eqno{(1)}$$ Every Moufang loop is
an $IP$-loop. The \textit{inner mappings} $T(a), R(a,b), L(a,b)$
are defined by
$$T(a) = L(a)^{-1}R(a), \quad R(a,b) = R(ab)^{-1}R(b)R(a),$$
$$L(a,b) = L(ab)^{-1}L(a)L(b), \eqno{(2)}$$
where $R(a)x = xa, L(a)x = ax$. The subloop $H$ of loop $Q$ is
called \textit{normal} in $Q$, if
$$xH = Hx,\quad x\cdot yH = xy\cdot H, \quad H\cdot xy = Hx\cdot y \eqno{(3)}$$
or by (2)
$$T(x)H = H,\quad L(x,y)H = H,\quad R(x,y)H = H \eqno{(4)}$$
for every $x, y \in Q$. The \textit{center} $Z(Q)$ of loop $Q$ is
a normal subloop $Z(Q) = \{x \in Q \vert x\cdot yz = xy\cdot z,
zy\cdot x = z\cdot yx, xy = yx$ $\forall y, z \in Q\}$. If $Z_1(Q)
= Z(Q)$, then the normal subloops $Z_{i+1}(Q): Z_{i+1}(Q)/Z_i(Q) =
Z(Q/Z_i(Q))$ are inductively determined. Loop $Q$ is called
\textit{centrally nilpotent of class $n$}, if its \textit{upper
central series} has the form $\{1\} \subset Z_1(Q) \subset \ldots
\subset Z_{n-1}(Q) \subset Z_n(Q) = Q$.

We pass to examining these loops. Let $Q$ be an arbitrary loop and
let $a, b, c \in Q$. The solution of the  equation  $ab\cdot c =
ax\cdot bc$ (respect. $c\cdot ba = cb\cdot xa$) is denoted by
$\alpha(a,b,c)$ (respect. $\beta(a,b,c)$) and is called the
\textit{associator of type $\alpha$} (respect.  \textit{of type
$\beta$}) of elements $a, b, c \in Q$. The \textit{commutator}
$(a,b)$ of elements $a, b \in Q$ is determined by the equality $ab
= b\cdot a(a,b)$. By (2) these definitions can be written in the
following way:
$$T(b)a = a(a,b),\quad R(b,c)a = a\alpha(a,b,c),\quad L(c,b)a =
\beta(a,b,c)a. \eqno{(5)}$$

\textbf{Lemma 1.} \textit{Let $N$ be a normal subloop of $IP$-loop
$Q$. Then the subloop $H$, generated by all elements of form
$\alpha(n,x,y), \beta(n,x,y), (n,x)$, where $n \in N, x, y \in Q$
is normal in $Q$ and $H \subseteq N$.}
\smallskip\\

\textbf{Proof.} From (5) we get $\alpha(n,x,y) = L(n)^{-1}R(x,y)n,
\beta(n,x,y) = \break R(n)^{-1}L(x,y), (n,x) = L(n)^{-1}T(x)n$.
Subloop $N$ is normal in $Q$, then by (4) $\alpha(n,x,y),
\beta(n,x,y), (n,x) \in N$ for any $n \in N$ and any $x, y \in Q$.
Then $H \subseteq N$. Let $h \in H$. By (5) $T(x)h = h(h,x) \in H,
L(x,y)h = \beta(h,y,x)h \in H, R(x,y)h = h\alpha(h,x,y) \in H$.
Then by (2) $ha \in aH$, $ha\cdot b \in H\cdot ab$, $a\cdot bh \in
ab\cdot H$. $Q$ is an $IP$-loop. Thus $a^{-1}h^{-1} \in Ha^{-1}$,
$b^{-1}\cdot a^{-1}h^{-1} \in b^{-1}a^{-1}\cdot H$,
$h^{-1}b^{-1}\cdot a^{-1} \in H\cdot b^{-1}a^{-1}$. Hence $aH =
Ha$, $a\cdot bH = ab\cdot H$, $Ha\cdot b = H\cdot ab$ for all $a,
b \in Q$ and by (3) $H$ is normal in $Q$.  This completes the
proof of Lemma 1.
\smallskip\\

Let $Q$ be an arbitrary loop. We note $Q_0 = Q$ and by induction
we determine $Q_{i+1}$. This is a normal subloop of loop $Q$,
generated by all the expressions $\alpha(n,x,y)$, $\beta(n,x,y)$,
$(n,x)$, where $n \in Q_i, x, y \in Q$. Then, again by induction
we obtain a series of normal subloops, by Lemma 1 $$Q = Q_0
\supseteq Q_1 \supseteq \ldots \supseteq Q_i \supseteq \ldots,$$
which we call the \textit{lower central series of loop $Q$}.

We observe that in [11] the associator $[a,b,c]$ and commutator
$[a,b]$ of elements $a, b, c \in Q$ are defined by the equalities
$ab\cdot c = (a\cdot bc)[a,b,c], ab = (ba)[a,b]$ for an  arbitrary
loop $Q$. Let now $Q$ be a Moufang loop and $a, b, c \in Q$. Loop
$Q$ is a $IP$-loop. Then $ab\cdot c = (a\cdot bc)[a,b,c]$,
$(ab\cdot c)[a,b,c]^{-1} = a\cdot bc$, $[a,b,c]^{-1} = (ab\cdot
c)^{-1}(a\cdot bc)$, $[a,b,c]^{-1} = (c^{-1}\cdot
b^{-1}a^{-1})(a\cdot bc)$ and by (1) we have $a[a,b,c]^{-1} =
a((c^{-1}\cdot b^{-1}a^{-1})(a\cdot bc)) = (a(c^{-1}\cdot
b^{-1}a^{-1})\cdot a)(bc) = ((ac^{-1})(b^{-1}a^{-1}\cdot a))(bc) =
(ac^{-1}\cdot b^{-1})(bc)$, i.e. $a[a,b,c]^{-1}\cdot c^{-1}b^{-1}
= ac^{-1}\cdot b^{-1}$. Further, $(a[a,b,c]^{-1}\cdot
c^{-1}b^{-1})^{-1} = (ac^{-1}\cdot b^{-1})^{-1}, bc\cdot
[a,b,c]a^{-1} = b\cdot ca^{-1}$. Then, from here and from the
definition of associator of type $\alpha$ or $\beta$, it follows
that $$[a,b,c]^{-1} = \alpha(a,b^{-1},c^{-1}),\quad [a,b,c] =
\beta(a^{-1},b,c). \eqno{(6)}$$ Further, the Moufang loop is
di-associative, so $[a,b] = (a,b)$. That is why from the
definition of subloop $Q_i$ and (7)  we get
\smallskip\\

\textbf{Prorosition 1.} \textit{The subloop $Q_{i+1}$ ($i =
0,1,\ldots$) of the lower central series of Moufang loop $Q$ is
generated by all the associators $[n,x,y]$ and all the commutators
$[n,x]$, where $n \in Q_i$, $x, y \in Q$.}
\smallskip\\

We will call the series of normal subloops $Q = C_0 \supseteq C_1
\supseteq \ldots \supseteq C_r = 1$ \textit{central}, if
$$C_i/C_{i+1} \subseteq Z(Q/C_{i+1}) \quad\text{for all}\quad i \eqno{(7)}$$
or, that is equivalent,
$$(C_i,Q)_{\alpha\beta} \subseteq C_{i+1} \quad \text{for all} \quad i, \eqno{(8)}$$
where $(C_i,Q)_{\alpha\beta}$ means a normal subloop of loop $Q$,
generated by all the elements of form $\alpha(u,x,y),
\beta(u,x,y), (u,x)$ ($u \in C_i, x, y \in Q$).
\smallskip\\

\textbf{Lemma 2.}  \textit{Let $Q = C_0 \supseteq C_1 \supseteq
\ldots \supseteq C_r = 1$ be a central series, $\{Z_i\}$ be the
upper central series, $\{Q_i\}$ be the lower central series of
loop $Q$. Then $C_{r - i} \subseteq Z_i$, $C_i \supseteq Q_i$, for
$i = 0, 1,\ldots,r$.}
\smallskip\\

\textbf{Proof.}  We have $C_0 = Q = Q_0$. Suppose that $C_i
\supseteq Q_i$. By (8) $(C_i,Q)_{\alpha\beta} \subseteq C_{i+1}$.
But then $Q_{i+1} = (Q_i,Q)_{\alpha\beta} \subseteq
(C_i,Q)_{\alpha\beta} \subseteq C_{i+1}$. We suppose now that
$C_{r - i} \subseteq Z_i$ for a certain $i$. Then loop $Q/Z_i$ is
the homomorphic image of loop $Q/C_{r - i}$ with kernel $Z_i/C_{r
- i}$. But by (7)
$$C_{r - i - 1}/C_{r - i} \subseteq Z(Q/C_{r - i}),$$
where from it follows that the homomorphic image of subloop $C_{r
- i - 1}/C_{r - i}$ must lie in the center $Z(Q/Z_i)$. It is clear
that this image is the subloop $(C_{r - i - 1} \cup Z_i)/Z_i$,
while $Z(Q/Z_i) = Z_{i + 1}/Z_i$. Consequently, $C_{r - i - 1}
\subseteq C_{r-i}\cup Z_i \subseteq Z_{i-1}$, as required.
\smallskip\\

\textbf{Prorosition 2.} \textit{Loop $Q$ is centrally nilpotent of
class $n$ if and only if its upper or lower central series have
respectively the forms}
$$1 = Z_0 \subset Z_1 \subset \ldots \subset Z_n = Q,\quad Q = Q_0
\supset Q_1 \supset \ldots \supset Q_n = 1.$$
\smallskip\\

\textbf{Proof.}  The statement of theorem for upper central series
follows from the definition of centrally nilpotent loop. Further,
if a central series of length $n$ exists, then from Lemma 2 it
follows that the lengths of the upper and lower central series do
not exceed $n$. But as there is a term by term inclusion between
the elements of these series, their  lengths  are  equal, and the
series have the indicated form, as required.

\section{Moufang loops and alternative algebras}

By analogy to Lemma 1 from [2] it is proved.
\smallskip\\

\textbf{Lemma 3.} \textit{Let $A$ be an alternative algebra and
let $Q$ be a subloop of $U(A)$. Then the restriction of any
homomorphism of algebra $A$  upon $Q$ will be a homomorphism on
the loop. Hence any ideal of $A$ induces a  normal subloop  of
$Q$.}
\smallskip\\

Let $L$ be a free Moufang loop, let $F$ be a field and let $FL$ be
a loop algebra of loop $L$ over field $F$. We remind that $FL$ is
a free module with basis $\{g \vert g \in L\}$ and the
multiplication of elements of the basis is defined by their
multiplication in loop $L$. Let $(u,v,w) = uv\cdot w - u\cdot vw$
denote the associator of elements $u, v, w$ of algebra $FL$. We
denote by $I$ the ideal of $FL$, generated by the set

$$\{(a,b,c) + (b,a,c), (a,b,c) + (a,c,b) \vert \forall a, b, c \in L\}.$$ It
is shown in [2] that algebra $FL/I$ is alternative and loop $L$ is
embedded (isomor\-phically) in the loop $U(A)$. Further we
identify the loop $L$ with its isomorphic image in $U(FL/I)$.
Without causing any misunderstandings, like in [2], we will denote
by $FL$ the quotient algebra $FL/I$ and call it $''$loop
algebra$''$ (in inverted commas).  Further, we will identify the
field $F$ with subalgebra $F1$ of algebra $FL$, where $1$ is the
unit of loop $L$.

Let now $Q$ be an arbitrary Moufang loop. Then $Q$ has a
representation as a quotient loop $L/R$ of the free Moufang loop
$L$ by the normal subloop $R$. Sums $\sum_{g\in L}\alpha_gg$, are
elements of algebra $FL$, where $\alpha_g \in F$. Let us determine
the homomorphism $\eta$ of alternative $F$-algebra $FL$ by the
rule: $\eta(\sum_{g \in L}\alpha_gg) = \break \sum_{g \in
L}\alpha_gRg$. Sums $\sum_{q\in Q}\alpha_qq$, are elements of
algebra $\eta(FL)$. We denote $\eta(FL) = FQ$ and the alternative
algebra $FQ$ will be called \textit{$''$loop algebra$''$} (in
inverted commas) of loop $Q$. Let now $H$ be a normal subloop of
$Q$ and let $\{u_i\}$ be a full  representative system of the
cosets  $Q/H$. Then any element $y$ in $FQ$ can be presented in
the form
$$y = x_1u_1 + \ldots + x_su_s, \eqno{(9)}$$ where $x_i = \sum_{h
\in H}\alpha_h^{(i)}h$.

We denote by $\omega H$ the ideal of $''$loop algebra$''$ $FQ$,
genera\-ted by the elements $1 - h$ ($h \in H$). Let $\varphi$ be
the homomorphism $Q \rightarrow Q/H$ of loops and we consider the
homomorphism $\overline\varphi$ of algebra $FQ$ defined by  the
rule: $\overline\varphi(\sum_{q \in Q}\alpha_qq) = \sum_{q \in
Q}\alpha_q\varphi(q)$. If $h \in H$ then $\overline\varphi(g(1 -
h)) = \varphi(g)(1 - \varphi(h)) = 0$, i.e.
$$\omega H \subseteq \text{Ker}\overline\varphi. \eqno{(10)}$$
Moreover, it is true.
\smallskip\\

\textbf{Lemma 4.} \textit{Let $H$ be a normal subloop of Moufang
loop $Q$. Then $h \in H$ if and only if $1 - h \in \omega H$.}
\smallskip\\

\textbf{Proof.}  If $q \notin H$, then $Hq \neq H$ and $\overline
\varphi(1 - q) = H - Hq \neq 0$. Hence $1 - q \notin
\text{Ker}\overline\varphi$. But by (10) $\omega H \subseteq
\text{Ker}\overline\varphi$. Then $1 - q \notin \omega H$. This
completes the proof of Lemma 4.
\smallskip\\

\textbf{Lemma 5.} \textit{Let $H$ be a proper normal subloop of
Moufang loop $Q$ and let $\overline\varphi$ be the homomorphism of
$''$loop algebra$''$ $FQ$ induced by $H$. Then the following
statements are equivalent:}

\textit{1) $\omega H \subset \text{Ker}\overline\varphi$;}

\textit{2) there exists an element $y \in
\text{Ker}\overline\varphi$ such that  in representation (9) there
exists  such an element $x_i = \sum_{h \in H}\alpha_h^{(i)}h$ that
$\sum_{h \in H}\alpha_h^{(i)} \neq 0$;}

\textit{3) $\text{Ker}\overline\varphi = FQ$.}
\smallskip\\

\textbf{Proof.}  $1) \Rightarrow 2)$. Let $y \in
\text{Ker}\overline\varphi$ and we suppose that in representation
(9) $\sum_{h \in H}\alpha_h^{(i)} = 0$ for $i = 1, \ldots, s$.
Then $x = \sum_{h \in H}\alpha_h^{(i)}h = -\sum_{h \in
H}\alpha_h^{(i)}(1 - h) + \sum_{h \in H}\alpha_h^{(i)} = -\sum_{h
\in H}\alpha_h^{(i)}(1 - h)$. Hence $x_i \in \omega H$, and $y \in
\omega H$ as well. But this contradicts the strict inclusion
$\omega H \subset \text{Ker}\overline\varphi$. Hence for some $i$
$\sum_{h \in H}\alpha_h^{(i)} \neq 0$.

$2) \Rightarrow 3)$. We denote by $\overline u$ the image of $u
\in FQ$ in $FQ/\text{Ker} \overline \varphi$. It is clearly that
if $\overline a = \overline 0$  for some $a \in Q$  then
$\text{Ker}\overline\varphi = FQ$. Let (9) be such a
representation of element $y$ of item 2) that the number $s$ of
representatives $u_i$ is minimal. We denote $\beta_i = \sum_{h \in
H}\alpha_h^{(i)}$ and let $i = 1$. Then $\overline 0 = \overline
\varphi(y) = \overline \varphi(x_1)\overline \varphi(u_1) + \ldots
+ \overline \varphi(x_s)\overline \varphi(u_s) = \sum_{h \in
H}\alpha_h^{(1)}\varphi(u_1) + \ldots + \sum_{h \in
H}\alpha_h^{(s)}\varphi(u_s)$, $\beta_1\varphi(u_1) =
-\beta_2\varphi(u_2) - \ldots - \beta_s\varphi(u_s)$,
$\varphi(u_1) = \gamma_2\varphi(u_2) + \ldots +
\gamma_s\varphi(u_s)$ ($\gamma_i = -\beta_i/\beta_1$), $Hu_1 =
\gamma_2Hu_2 + \ldots + \gamma_sHh_s$, $u_1 = \gamma_2h_2u_2 +
\ldots + \gamma_sh_su_s$. We substitute in (9) the expression
obtained for $h_1$. If $s > 1$ we get that the element $y$ has a
representation of type (9) with less  representatives that $s$.
But is contradicts the minimum of number  $s$. Hence $y =
\beta_1u_1$. But $\overline y = \overline 0$. Then   $\text{Ker}
\overline \varphi = FQ$.

$3) \Rightarrow 1)$ follows from Lemma 4 as the subloop $H$ is
proper. This completes the proof of Lemma 5.
\smallskip\\

\textbf{Lemma 6.} \textit{Let $H$ be a proper normal subloop of
free loop $L$, let $\overline\varphi$ be the homomorphism of
$''$loop algebra$''$ $FL$ induced by $H$ and let $\omega H \subset
\text{Ker}\overline\varphi$. If $I_1, I_2$ are a proper ideals of
algebra $FL/\omega H$ then and $I_1 + I_2$ is also proper ideal of
$FL/\omega H$.}
\smallskip\\

\textbf{Proof.}  As $\omega H \subset \text{Ker}\overline\varphi$
then $\omega H$ is a proper ideal of $FL$ and by Lemma 5
$FL/\omega H$ is a non-trivial algebra. By Lemma 3 $\omega H$
induces a normal subloop $R \subset H$ of loop $L$ and $F(L/R) =
FL/\omega H$. Let $L/R = Q$. Again by Lemma 3 the ideals $I_1$ and
$I_2$ induce  a normal subloops $K_1$ and $K_2$ of $Q$
respectively. $I_1$ and $I_2$ are proper ideals of $FQ$, then by
Lemma 5 $\omega K_1 = I_1$ and $\omega K_2 = I_2$. Any element $x$
in $FQ$ has the form $x = \sum\alpha_iq_i$, where $\alpha_i \in
F$, $q_i \in Q$. If $x \in \omega K_1$, then from definition of
ideal $\omega K_1$ it follows that $\sum\alpha_i = 0$.
Analogically, if $y \in \omega K_2$ and $y = \sum\beta_jq_j$ then
$\sum\beta_j = 0$. Hence any element $z$ in $\omega K_1 + \omega
K_2$ has a form $z = \sum\gamma_kq_k$ with $\sum\gamma_k = 0$.
Hence $\omega K_1 + \omega K_2 \neq FQ$, thus $I_1 + I_2$ is a
proper ideal of $FQ$. This completes the proof of Lemma 6.
\smallskip\\

\textbf{Lemma 7.} \textit{Let $Q$ be a Moufang loop, let $Q =
L/H$, where $L$ is a free Moufang loop and we suppose that $FQ =
FL/\omega H$. Then the loop $Q$ can be embedded into a loop of
invertible elements $U(FQ)$.}
\smallskip\\

\textbf{Proof.}  In accordance with [2] we consider that $L$ is a
subloop of $U(FL)$ of $''$loop algebra$''$ $FL$. By Lemma 3 the
ideal $\omega H$ of $FL$ induces a normal subloop $R$ of $U(FL)$
and $F(L/R) = FL/\omega H$. Hence $L/R$ is a subloop of
$U(FL/\omega H)$. Let $\overline \varphi$ be a homomorphism $FL
\rightarrow F(L/R)$ induced by homomorphism $L \rightarrow L/R$.
By (10) we have $\omega R \subseteq \text{Ker}\overline\varphi$.
But $\text{Ker}\overline\varphi = \omega H$ and $FL/\omega H$ is a
non-trivial algebra. Then by Lemma 5 $\omega R =
\text{Ker}\overline\varphi$ and by Lemma 4 from $\omega R = \omega
H$ it follows that $R = H$. This completes the proof of Lemma 7.

Let $F$ be a field. Let us consider a classical matrix
Cayley-Dickson algebra $C(F)$. It consists of matrices of form
$$a = \left(
\begin{array}{ll} \alpha_1 & \alpha_{11}\\
\alpha_{21} & \alpha_2  \end{array} \right), \eqno{(11)}$$ where
$\alpha_1, \alpha_2 \in F$, $\alpha_{12}, \alpha_{21} \in F^3$.
The addition and multiplication by scalar of elements of algebra
$C(F)$ is represented by ordinary addition and multiplication by
scalar of matrices, and the multiplication of elements of algebra
$C(F)$ is defined by the rule

$$ \left(
\begin{array}{ll} \alpha_1 & \alpha_{12}\\
\alpha_{21} & \alpha_2  \end{array} \right) \left(
\begin{array}{ll} \beta_1 & \beta_{12} \\ \beta_{21} & \beta_2
\end{array} \right) = $$
$$ \left( \begin{array}{ll} \alpha_1 \beta_1 + (\alpha_{12},
\beta_{21}) & \alpha_1 \beta_{12} + \beta_2 \alpha_{12} -
\alpha_{21} \times \beta_{21} \\ \beta_1 \alpha_{21} + \alpha_2
\beta_{21} + \alpha_{12} \times \beta_{12} & \alpha_2 \beta_2 +
(\alpha_{21}, \beta_{12}) \end{array} \right), $$ where for
vectors $\gamma = (\gamma_1, \gamma_2, \gamma_3),$ $\delta =
(\delta_1, \delta_2, \delta_3) \in A^3$  $(\gamma, \delta) =
\gamma_1\delta_1 + \gamma_2\delta_2 + \gamma_3\delta_3$ denotes
their scalar product and $\gamma \times \delta = (\gamma_2\delta_3
- \gamma_3\delta_2, \gamma_3\delta_1 - \gamma_1\delta_3,
\gamma_1\delta_2 - \gamma_2\delta_1)$ denotes the vector product.
Algebra  $C(F)$ is alternative. It is also split and quadratic
over $F$, i.e. each element $a \in C(F)$ satisfies the identity

$$a^2 - t(a)a + n(a) = 0, n(a), t(a) \in F$$ and admits
composition, i. e.

$$n(ab) = n(a)n(b)$$ for $a, b \in C(F)$. Track $t(a)$ and norm
$n(a)$ are defined by the equalities $t(a) = \alpha_1 + \alpha_2,$
$n(a) = \alpha_1\alpha_2 - (\alpha_{12}, \alpha_{21})$.

We denote $M_0(F) = \{u \in C(F) \vert n(u) = 1\}$. It follows
from the relation $n(ab) = n(a)n(b)$ that if $a, b \in M_0(F)$
that $ab \in M_0(F)$. Further, for $a \in M_0(F)$  $a^2 - t(a)a +
1 = 0, -a^2 + t(a)a = 1, a(-a +t(a)) = 1$, i. e. $a$ has an
inverse element in $M_0(F)$. Therefore  $M_0(F)$ is a loop. We
denote $U(C(F)) = U(F)$. Analogically it is proved that $U(F) =
\{a \in C(F) \vert n(a) \neq 0\}$. Moufang identities hold in
alternative algebras, hence $U(F)$ is a Moufang loop. If $u, v \in
U(F)$ then from the relation $n(uv) = n(u)n(v)$ it follows that
the mapping $u \rightarrow n(u)$ is a homomorphism of the loop
$U(F)$ upon $F$. The inverse image of $1 \in F$ is $M_0(F)$. Hence
$M_0(F)$ is a normal subloop of $U(F)$. Let $Z(U)$ denote the
center of $U(F)$ and $Z(M_0)$ denote the center of $M_0(F)$. In
[6] is proved that $Z(M_0) = U(F) \cap Z(U)$, $Z(M_0)$ is
generated by element $-1$ and $Z(F)$ be made of all matrices of
form (11) for which $a_1 = a_2 \neq 0$, $\alpha_{12} = \alpha_{21}
= 0$. Then to within an isomorphism $M_0(F)/Z(M_0)$ is a normal
subloop of $U(F)/Z(U)$, $M_0(F)/Z(M_0) \cong U(F)/Z(U)$ if and
only if $U(F) = M_0(F)Z(U)$ and this will be true if the field $F$
is closed under the square root operation. This means that the
equation $x^2 - a = 0$ is solved in $F$ for all $0 \neq a \in F$.
Obviously, this equation is solved in the field of real numbers,
in the field of complex numbers and is unsolved in any simple
field, i.e., in the field of rational numbers and in the finite
field $FG(p)$, $p \neq 2$. We also mention that in [6] there is
constructed a Cayley-Dickson division algebra over the field of
all formal power series $\sum^{\infty}_{k = n}a_kt^k$ with real
coefficients and $n$ is either positive, negative, or zero, for
which the equation $x^2 - a = 0$ is unsolved. Let now $F$ be an
arbitrary field and $a \in F$. If $a$ is not a square in $F$ then
the polynomial $x^2 - a$ doesn't have a square in $F$, hence it is
irreducible. We suppose that $\text{char} F \neq 2$. Then the
polynomial $x^2 - a$ is separable as $a \neq 0$ and if $\alpha$ is
its square, then the extension $F(\alpha)$ is Galois. Its Galois
group $G$ is cyclic of order 2. The order of $G$ coincides with
the degree of extension $\vert F(\alpha) : F \vert$. The numbers
$\pm 1$ are square roots of the unity element and they belong only
to $F$. Then it follows from [14, pag. 216] that for $\text{char}
F \neq 2$ the square extension of field $F$ is equivalent with
connection of the square root of some element of $F$.
Consequently, the field $F$ is closed under square root operation
if and only if $\sqrt{a} \in F$ for some $a \in F$.

Let now $FG(p^n)$ be a finite field, where $p \neq 2$. We consider
the simple field $FG(p)$ as basic for $FG(p^n)$. Let  $G$ be a
Galois group of extension $FG(P^n)$ over $FG(p)$. Then $G$ will be
a cyclic group of order $n$. Let $n = m_1\ldots m_k$ be a
decomposition of $n$ in prime factors. Then $G$ has such a
composition series $1 = G_0 \subset G_1 \subset \ldots \subset
G_{m_k} = G$ that $G_i/G_{i-1}$, $i - 1, \ldots, m_k$, is a cyclic
group of order $m_i$. By the main theorem of Galois theory chain
of subfields of field $FG(p^n)$ corresponds to it, where each next
item has the degree $m_i$ over the preceding one. It follows from
the aforementioned reasonings that \textit{a finite field
$FG(p^n)$, $p \neq 2$, is closed under square root operation if
and only if the exponent $n$ is an even number.}

We denote $M(F) = M_0(F)/Z(M_0)$. It takes place.
\smallskip\\

\textbf{Lemma 8.} [7, 8]. \textit{Let $P$ be an algebraically
closed field. Only and only the loops $M(F)$ of the matrix
Cayley-Dickson algebra $C(F)$, where $F$ is a subfield of field
$P$ are with precise till isomorphism non-associative simple
Moufang loops.}
\smallskip\\

\textbf{Theorem 1.} \textit{Let $P$ be an algebraically closed
field and let $F$ be a subfield of $P$. The simple loops $M(F)$
they and only they are not embedded into a loop of type $U(A)$ of
any unitaly alternative algebras $A$ if $\text{char} F \neq 2$ and
$F$ is closed under square root operation. The remaining loops can
be embedded into a loop of type $U(A)$. In particular, any
non-simple loop $Q$ can be embedded into a loop $U(FQ)$, where $F$
is an arbitrary field and $FQ$ is the $''$loop algebra$''$ of
$Q$.}
\smallskip\\

\textbf{Proof.}  Let $Q$ be an arbitrary Moufang loop. $Q$ has a
representation $Q = L/H$, where $L$ is a free Moufang loop. Let
$F$ be an arbitrary field and let $P$ be the algebraic closing of
$F$. We consider the inclusion (10). It can be divided in two
cases: a) $\omega H = \text{Ker}\overline\varphi$; b) $\omega H
\subset \text{Ker}\overline\varphi$. If the case a) holds then by
Lemma 7 $Q$ can be embedded into loop $U(FQ)$ of $''$loop
algebra$''$ $FQ$.

Now we suppose the case b) holds. By Lemma 5 $\omega H$ is a
proper ideal of $FQ$ and by Lemma 3 $\omega H$ induces a normal
subloop $R \subset H$ of $U(PQ)$. We denote $\overline Q = Q/R$
and let $S$ be the ideal of $P\overline Q$,  generated by all
proper ideals $J_i$ of algebra $P\overline Q$. Let us show that
$S$ is also proper ideal of algebra $P\overline Q$. Indeed, ideal
$S$ is the subgroup of additive group of algebra $P\overline Q$,
consisting of all possible finite sums  $\sum \alpha_ju_j$, where
$\alpha _j \in P, u_j \in J_i$. Let us suppose that for any
elements $x \in \overline Q$ in $P\overline Q$ there is such a
finite number of ideals $J_i$, that $x \in \sum J_i$. The algebra
$P\overline Q$ is generated as a $P$-module by elements $x \in
\overline Q$. Then $\sum J_i = P\overline Q$. But this contradicts
Lemma 6. Therefore this case is impossible.

Let us now consider the second possible case. Let there exist be
such ideals $J_1, \ldots, J_k$ that for element $1 \neq a \in
\overline Q$ $a \in \sum J_i$ and let us suppose that for element
$b \in \overline Q$ $b \notin \sum J_i$. We denote by $T$ the set
of all ideals of $P\overline Q$, containing the element $a$, but
not containing the element $b$. By Zorn's Lemma there is a maximal
ideal $I_1$ in $T$. We denote by $I_2$ the ideal of algebra
$P\overline Q$, generated by all proper ideals of $P\overline  Q$,
that don't belong to ideal $I_1$. Then $S = I_1 + I_2$ and by
Lemma 6 $S$ is a proper ideal of $P\overline  Q$.

By Lemma 3 $S$ induces a normal sublooop $H$ of loop $\overline
Q$. We denote $\overline{\overline Q} = \overline Q/H$. The ideal
$S$ is proper, then by Lemma 5 $P\overline{\overline Q} =
P\overline Q/\omega H$ and the loop $\overline{\overline Q}$ is
embedded into a loop $U(P\overline{\overline Q})$. The alternative
algebra $P\overline{\overline Q}$ is simple and $P$ is an
algebraically closed field. Then by Kleinfeld Theorem it is a
matrix Cayley-Dickson algebra over its center. In [7, 8] is proved
that the loop $\overline{\overline Q}$ is isomorphic with a loop
$M_0(F)$ for an appropriate subfield $F$ of  $P$ and the loop $Q$
is isomorphic with loop $M(F) = M_0(F)/<-1>$. From Lemma 8 it
follows that if case b) holds thus the loop $Q$ is simple and if
the case a) holds thus the loop $Q$ is non-simple.

We consider the case b). If $\text{char}F = 2$ then the loop $Q
\cong M(F)$ is embedded into a loop $U(F\overline{\overline Q})$.
For $\text{char} F \neq 2$ we consider the subcases: c) $F$ is not
closed under square root operation; d) $F$  is closed under square
root operation. If the case c) holds then as is indicated before
Lemma 8 to within an isomorphism $Q$ is a proper normal subloop of
$U(F)/Z(U)$. Hence $U(F)/Z(U)$ is a non-simple loop and by case a)
$U(F)/Z(U)$ can be embedded into a loop of type $U(A)$. Then and
the loop $Q$ can be embedded into a loop of type $U(A)$. If the
case d) holds then $\sqrt{2} \in F$ and repeating almost word dy
word the proof of Theorem 1 from [3] it is proved that the simple
Moufang loop $Q \cong M(F)$  is not embedded into a loop of type
$U(A)$ for any unital alternative algebra $A$. This completes the
proof of Theorem 1.
\smallskip\\

\textbf{Corollary 1.} \textit{In the class of finite Moufang loops
the simple Moufang loops $M(F)$, where $F = FG(p^n)$, $p \neq 2$,
$n = 2k, k = 0, 1, 2, \ldots$, they and only they are not embedded
onto a loop of invertible elements of any unitaly alternative
algebra.}
\smallskip\\

\textbf{Proof.}  This statement follows from the remark made
before Lemma 8 and Theorem 1.

Let $Q$ be an arbitrary non-simple Moufang loop with set of
generators $\{g_1, g_2, \ldots \break \ldots, g_i \ldots\}$ and
let $L$ be the free Moufang loop with set of free generators
$\{x_1, x_2, \ldots \break \ldots, x_i, \ldots\}$. Then $Q = L/H$.
Temporarily, by $FL$, $FQ$ we denote the loop algebras and by
$\overline{FL}$, $\overline{FQ}$ we denote the $''$loop
algebras$''$ of loop $L$, $Q$ respectively. Let $I$ be the ideal
of $FL$ generated by set $\{(a, b, c) + (b, a, c), (a, b, c) + (a,
c, b) \vert \forall a, b, c \in L\}$ and let $J$ be the ideal of
$FQ$ generated by set $\{(u, v, w) + (v, u, w), (u, v, w) + (u, w,
v) \vert \forall u, v, w \in Q\}$. The $''$loop algebra$''$
$\overline{FL}$ we defined as $\overline{FL} = FL/I$, not the
$''$loop algebra$''$ $\overline{FQ}$ we defined differently,  as
$\overline{FQ} = FQ/\omega H$. However it take place.
\smallskip\\

\textbf{Prorosition 3.} \textit{Let $Q$ be a non-simple Moufang
loop. Then $\overline{FQ} = FQ/J$.}
\smallskip\\

\textbf{Proof.}  We suppose that a loop $K$ is embedded into a
loop of type $U(A)$ for some algebra $A$. Below the isomorphic
image of $K$ in $U(A)$ we identify with $K$. Then, according to
the Theorem 1 $Q \subseteq U(\overline{FQ})$ and any element in
$\overline{FQ}$ has a form of finite sum $\sum_{q\in Q}\alpha_qq$.
By definition, the algebra $FQ$ is a free $F$-module with basis
$\{q \vert q\in Q\}$. Then the mappings $q \rightarrow q$ induce a
homomorphism $\mu: FQ \rightarrow \overline{FQ}$. $\overline{FQ}$
is an alternative algebra, then $J \subseteq \text{Ker}\mu$. Thus,
the homomorphism $\mu$ induces the homomorphism $\varphi: FQ/J
\rightarrow \overline{FQ}$.

Further we denote by $\overline g_i$ the image of generators $g_i$
of $Q$ under the homomor\-phism $FQ \rightarrow FQ/J$. Let $\eta:
FL \rightarrow FQ/J$ be the homomorphism defined by mappings $x_i
\rightarrow \overline g_i$. It follows from definition of loop
algebra that $FQ/J$ is an alternative algebra. Thus $I \subseteq
\text{Ker}\eta$ and $\eta$ induces homomorphism $\xi: FL/I
\rightarrow FQ/J$. The loop $L$ is non-simple. Then by Theorem 1
$L \subseteq U(\overline{FL})$. We have  also shown  that $Q
\subseteq U(\overline{FQ})$. We have $Q = L/H$. Thus $1 - H
\subseteq \text{Ker}\xi$ or $\omega H \subseteq \text{Ker}\xi$.
Then the homomorphism $\xi$ induces a homomorphism $\psi:
\overline{FL}/\omega H \rightarrow FQ/J$. Not
$\overline{FL}/\omega H = \overline{FQ}$. Hence we have a
homomorphism $\psi: \overline{FQ} \rightarrow FQ/J$ which together
with homomorphism $\varphi: FQ/J \rightarrow \overline{FQ}$ show
that $\varphi, \psi$ are an isomorphisms. This completes the proof
of Proposition 3.
\smallskip\\

Let $Q$ be a loop Moufang, let $F$ be an arbitrary field and let
$FQ$ be the $''$loop algebra$''$. As in [2] the ideal of $FQ$
generated by set $\{1 - q \vert \forall q \in Q\}$ will be called
\textit{$''$augmentation ideal$''$} of $FQ$ and will be denoted by
$\omega Q$. It is easily to see that $\omega Q = \text{Ker}
\varphi$, where $\varphi$ is the homomorphism of $FQ$ upon $F$
determined by rule $\varphi(\sum_{q \in Q}\alpha_qq) = \sum_{q \in
Q}\varphi_q$.
\smallskip\\

\textbf{Prorosition 4.} \textit{Let $Q$ be a non-simple Moufang
loop, let $H, H_1, H_2$ be its normal subloops, let $F$ be an
arbitrary field and let $FQ$ and $\omega Q$ are respectively
''loop algebra'' and ''augmentation ideal'' of $Q$. Then}

\textit{1) $FQ$ is generated as $F$-module by set $\{q \in Q \vert
\forall q \in Q\}$,}

\textit{2) $\omega Q$ is generated as $F$-module by set $\{1 - q
\vert \forall q \in Q\}$,}

\textit{3) $\omega Q = \{\sum_{q \in Q}\lambda_qq \vert \sum_{q
\in Q}\lambda_q
 = 0\}$,}

\textit{4) $FQ/\omega H \cong F(Q/H), \omega Q/\omega H \cong
\omega (Q/H)$,}

\textit{5) if the elements $h_i$ generate the subloop $H$, then
the elements $1 - h_i$ generate the ideal $\omega H$; if $H_1 \neq
H_2$, then $\omega H_1 \neq \omega H_2$; if $H_1 \subset H_2$,
then $\omega H_1 \subset \omega H_2$; if $H = \{H_1,H_2\}$, then
$\omega H = \omega H_1 + \omega H_2$,}

\textit{6) $F \cap \omega\overline{Q} = 0$, $F\overline{Q} = F +
\omega\overline{Q}.$}
\smallskip\\

\textbf{Proof.}  1). It follows from Theorem 1 and Proposition 3.

2). As $(1 - q)q^{\prime}= (1 - qq^{\prime}) - (1 - q^{\prime})$,
then the $''$augmentation ideal$''$ $\omega Q$ is generated by the
elements of form $1 - q$, where $q \in Q$.

3). Denote $R = \{\sum_{q \in Q}\lambda_qq \vert \sum_{q \in
 Q}\lambda_q = 0\}$. Obviously, $\omega Q \subseteq R$. Conversely,
 if $r \in R$ and $r = \sum_{q \in Q}\lambda_qq$, then $-r =
- \sum_{q \in Q}\lambda_qq = (\sum_{q \in Q}\lambda_q)1 - \sum_{q
\in Q}\lambda_qq = \sum_{q \in Q}\lambda_q(1-q) \in \omega Q$,
i.e. $R \subseteq \omega Q$. Hence $R = \omega Q$.

4). The relation $F(Q/H) \cong FQ/\omega H$ is the case a) of
proof of Theorem 1. The mapping $\varphi:FQ \rightarrow
F(Q)/\omega H$ keeps the sum of coefficients then by 3) from the
first  relation follows the second  relation of 4).

5). Let elements $\{h_i\}$ generate subloop $H$ and $I$ be an
ideal, generated by the elements $\{1 - h_i\}$. Obviously $I
\subseteq \omega H$. Conversely, let $g \in H$ and $g = g_1g_2$,
where $g_1, g_2$ are words from $h_i$. We suppose that $1 - g_1, 1
- g_2 \in I$. Then $1 - g = (1-g_1)g_2 + 1 - g_2 \in I$,  i.e. $I
= \omega H$. Let $H_1 \neq H_2$ and $g \in H_1, g \notin H_2$.
Then by  Lemma 4 $1 - g \in \omega H_1$, but $1 - g \notin \omega
H_2$. If $H = \{H_1, H_2\}$, then by the first statement of 5)
$\omega H = \omega H_1 + \omega H_2$.

6). From case a) of proof of Theorem 1 and Lemma 5 it follows that
$\omega Q \neq FQ$. We define the homomorphism of $F$-algebras
$\varphi$: $F\overline Q \rightarrow F$ by the rule
$\varphi(\sum\alpha_qq) = \sum\alpha_q$. We have that $\text{Ker}
\varphi = \omega Q$. Then $FQ = F + \omega Q$ and by 3) $\omega Q
\cap F = 0$. This completes the proof of Proposition 4.

\section{Finite Moufang $p$-loops}

\textbf{Lemma 9.} [11]. \textit{If $Q$ is a simple Moufang loop,
the order of every element of $Q$ devices the order of $Q$.}
\smallskip\\

\textbf{Lemma 10.}  \textit{Any finite Moufang $p$-loop $Q$ can be
embedded into a loop $U(FQ)$ of $''$loop algebra$''$ $FQ$.}
\smallskip\\

\textbf{Proof.}  If $Q$ is a non-simple loop then by Theorem 1 $Q$
can be embedded into a loop $U(FQ)$. Let now $Q$ be a simple
$p$-loop. By Lemma 8 $Q$ is isomorphic with loop $M(F) =
M_0(F)/<-1>$. $<-1>$ is cyclic group of order 2. Let $1, a_1,
\ldots, a_7$ the the canonical basis for $C(F)$. Then $a_i \in
M_0(F)$ and as $a_i^2 = -1$ thus $M(F)$ contains  $2$-elements. By
Lemma 9 $M(F)$ is $2$-loop. Then $\text{char} F = 2$ and by
Theorem 1 $Q$ can be embedded into a loop $U(A)$, where $A$ is an
algebra of type $FQ$. This completes the proof of Lemma 10.
\smallskip\\

If $A$ is a arbitrary $F$-algebra, then its $n$ degree $A^n$ is
$F$-module with a basis, consisting of products of any its $n$
elements with  any brackets distribution. Algebra $A$ is called
\textit{nilpotent} if $A^n = (0)$ for a certain $n$.
\smallskip\\

\textbf{Lemma 11.} \textit{Let $Q$ be a finite Moufang $p$-loop
and $F$ be a field  of characteristic $p$. Then the ''augmentation
ideal'' $\omega Q$ of ''loop algebra'' $FQ$ is nilpotent.}
\smallskip\\

\textbf{Proof.}  In accordance with Lemma 10 we consider that $Q
\subseteq FQ$ and we used 1) of Proposition 4. If $g \in Q$, then
by Lemma 9 $g^k = 1$, where $k = p^n$. We have $(1 - g)^k = 1 -
{1\choose k}g + \ldots + (-1)^i{i\choose k}g^i + ... + (-1)^kg^k$.
All binomial coefficients ${i\choose k}$ can be divided by $p$,
therefore $(1 - g)^k = 1 + (-1)^kg^k$. If $p = 2$, then $(1 - g
)^k = 1 + g^k = 0$, because $F$ is a field  of characteristic 2.
But if $p > 2$, then $(1 - g)^k = 1 - g^k = 0$. Then to the
algebra $\omega Q$ one can apply the statement: any alternative
$F$-algebra, generated as $F$-module by a finite set of nilpotent
elements, is nilpotent [18, pages 144, 408]. Consequently the
''augmentation  ideal'' $\omega Q$ is nilpotent, as required.
\smallskip\\

Let now $A$ be an alternative $F$-algebra with  unit $1$ and $B$
be a subalgebra from $A$, satisfying the law
$$x^m = 0. \eqno{(12)}$$
Then $1 - B = \{1 - b \vert b \in B\}$ will be a loop and $(1 -
b)^{-1} = 1 + b + \ldots + b^{m-1}$. We note $\sum x = 1 + x +
\ldots + x^{m-1}$. We remind that inscription $(a,b,c) = ab\cdot c
- a\cdot bc, (a,b) = ab - ba$ mean the associator and commutator
in algebra, but $[a,b,c] = (a\cdot bc)^{-1}\cdot (ab\cdot c),
[a,b] = a^{-1}b^{-1}\cdot ab$ are associator and commutator in
$IP$-loop.
\smallskip\\

\textbf{Lemma 12.} \textit{Let $A$ be an alternative algebra with
unit $1$ and $B$ its subalgebra, satisfying the law (18). Then for
$u, v, w \in B$ $[1 - u , 1 - v, 1 - w] = 1 - ((1 + w + \ldots +
w^{m-1})(1 + v + \ldots + v^{m-1})\cdot(1 + u + \ldots +
u^{m-1}))(u,v,w), [u,v] = (1 + u + \ldots
 + u^{m-1})(1 + v + \ldots + v^{m-1})(u,v)$.}
\smallskip\\

\textbf{Proof.}  We denote $1 - u = a, 1 - v = b, 1 - w = c$. Then
we have $[1 - u, 1 - v, 1 - w] = (a\cdot bc)^{-1}(ab\cdot c) =
(a\cdot bc)^{-1}(ab\cdot c) - (a\cdot bc)^{-1}(a\cdot bc) + 1 = 1
+ (a\cdot bc)^{-1}(a,b,c) = 1 + (((1 - w)^{-1}\cdot (1 -
v)^{-1})(1 - u)^{-1})(1 - u, 1 - v, 1 - w) = 1 - ((1 - w)^{-1}(1 -
v)^{-1}\cdot (1 - u)^{-1})(u,v,w) = 1 - ((1 + w + \ldots +
w^{m-1})(1 + v + \ldots + v^{m-1})\cdot (1 + u + \ldots +
u^{m-1}))(u,v,w)$.  The  second  equality  is  proved  by analogy.
\smallskip\\

\textbf{Lemma 13.} \textit{Let $Q$ be a Moufang loop and let the
''augmentation ideal'' $\omega Q$ be nilpotent. Then loop $Q$ is
centrally nilpotent.}
\smallskip\\

\textbf{Proof.}  It follows from  the definition of the
''augmentation ideal'' that $Q = Q_0 \subseteq 1 - \omega Q$. We
suppose that $Q_{i-1} \subseteq 1 - (\omega Q)^i$. Then it follows
from the Proposition 1 and Lemma 12 that $Q_i \subseteq 1 -
(\omega Q)^{i+1}$. Algebra $\omega Q$ is nilpotent and we suppose
that $(\omega Q)^{k+1} = (0)$. Then $Q_k = 1$ and by Theorem 1 and
Proposition 1 loop $Q$ is centrally nilpotent, as required.
\smallskip\\

It follows from Lemmas 11 and 13.
\smallskip\\

\textbf{Prorosition 5.} \textit{Any finite Moufang $p$-loop is
centrally nilpotent.}
\smallskip\\
In the end, I thank prof. I. P. Shestakov who kindly offered his
manuscript before its release.

\smallskip
Tiraspol State University of Moldova$$
$$

$$ $$

The author's address:

Sandu Nicolae Ion

Deleanu str 1

Apartment 60

Kishinev MD-2071, Moldova

E-mail: sandumn@yahoo.com
\end{document}